\magnification\magstep1
\baselineskip = 18pt
\def \Bbb {\bf}
\def\square{\vcenter{\hrule height1pt
\hbox{\vrule width1pt height4pt \kern4pt
\vrule width1pt}
\hrule height1pt}}
\centerline{\bf Shadows of Convex bodies}\bigskip
\centerline{by}\bigskip
\centerline{Keith Ball$^{(1)}$}
\centerline{Trinity College}
\centerline{Cambridge}\smallskip
\centerline{and}\smallskip
\centerline{Texas A\&M University}
\centerline{College Station, Texas}\bigskip

\noindent {\bf Abstract.} It is proved that if $C$ is a convex body in
${\Bbb R}^n$ then $C$ has an affine image $\widetilde C$ (of non-zero
volume) so that if $P$ is any 1-codimensional orthogonal projection,

$$|P\widetilde C| \ge |\widetilde C|^{n-1\over n}.$$

\noindent It is also shown that there is a pathological body, $K$, all of
whose orthogonal projections have volume about $\sqrt{n}$ times as large as
$|K|^{n-1\over n}$.\vskip3in

\noindent A.M.S. (1980) Subject Classification: \ 52A20, 52A40

\noindent  $^{(1)}$Supported in part by N.S.F. DMS-8807243

\vfill\eject

\noindent {\bf \S 0. Introduction.}

The problems discussed in this paper concern the areas of shadows
(orthogonal projections) of convex bodies and, to a lesser extent, the
surface areas of such bodies. If $C$ is a convex body in ${\Bbb R}^n$ and
$\theta$ a unit vector, $P_\theta C$ will denote the orthogonal projection
of $C$ onto the 1-codimensional space perpendicular to $\theta$. Volumes
and areas of convex bodies and their surfaces will be denoted with $|\cdot|$.

The relationship between shadows and surface areas of convex bodies is
expressed in Cauchy's well-known formula. For each $n\in {\Bbb N}$, let
$v_n$ be the volume of the $n$-dimensional Euclidean unit ball and let
$\sigma = \sigma_{n-1}$ be the rotationally invariant probability on the
unit sphere $S^{n-1}$. Cauchy's formula states that if $C$ is a convex body
in ${\Bbb R}^n$ then its surface area is

$$|\partial C| = {nv_n\over v_{n-1}} \int_{S^{n-1}} |P_\theta C|
d\sigma(\theta).$$

The classical isoperimetric inequality in ${\Bbb R}^n$ states that any body
has surface area at least as large as an Euclidean ball of the same volume.
The first section of this paper is devoted to the proof of a ``local''
isoperimetric inequality showing that all bodies have large shadows rather
than merely large surface area (or average shadow). The principal
motivation for this result is its relationship to a conjecture of Vaaler
and the important problems surrounding it. This theorem
 and its connection with Vaaler's conjecture are described at the
beginning of Section 1.

An important role is played in the theory of convex bodies by the so-called
``projection body'' of a convex body. It is easily seen, by considering
polytopes, that for every convex body $C\subset {\Bbb R}^n$, there is a
Borel measure $\mu$ on $S^{n-1}$ so that for each $\theta \in S^{n-1}$,

$$|P_\theta C| = \int _{S^{n-1}} |\langle \theta, \phi\rangle |
d\mu(\phi).$$

\noindent Hence, there is a norm $\|\cdot \|$ on ${\Bbb R}^n$, with
$\|\theta\| = |P_\theta C|$ for each $\theta \in S^{n-1}$, with respect to
which ${\Bbb R}^n$ is isometrically isomorphic to a subspace of $L_1$. The
unit ball of this norm is a symmetric convex body which will be denoted
$\Pi^*(C)$. The map $\Pi^*$, from the collection of  convex bodies in
${\Bbb R}^n$ to the collection of unit balls of representations of
$n$-dimensional subspaces of $L_1$ on ${\Bbb R}^n$, has been extensively
studied: \ see e.g. [B-L]. The restriction of $\Pi^*$ to the class of
centrally symmetric convex bodies was shown to be injective by Aleksandrov:
\  if $C$ and $D$ are centrally symmetric convex bodies and $|P_\theta C| =
|P_\theta D|$ for all $\theta \in S^{n-1}$ then $C = D$. If the condition
of central symmetry is dropped, $C$ and $D$ may not even be congruent: \
the Rouleaux triangle in ${\Bbb R}^2$ has all 1-dimensional shadows equal
in length to those of some disc.

The map $\Pi^*$ was shown to be surjective by Minkowski. What Minkowski's
proof actually gives (at least in the context of polytopes) is the result
stated in Section 2 as Lemma 6. (This rather detailed statement of
Minkowski's theorem will be needed for the construction of a pathological
body with large shadows.)

An important observation of Petty [P], is that if $T$ is a linear
operator on ${\Bbb R}^n$ of determinant 1 then, for every $C$,

$$ \Pi^*(TC) =
T(\Pi^*(C)).\eqno (1)$$

\noindent Motivated in part by Aleksandrov's theorem on the injectivity of
$\Pi^*$, Shephard asked whether, if $C$ and $D$ are centrally symmetric
convex bodies with

$$|P_\theta C| \ge |P_\theta D|  \quad \hbox{for all} \quad \theta \in
S^{n-1}$$

\noindent then necessarily $|C| \ge |D|$. This question was answered in the
negative by Petty and Schneider independently in [P] and [S]. (The
corresponding question for sections rather than shadows was posed in [B-P]
and answered (again in the negative) by [L-R].) The second section of this
paper contains a strongly negative answer to Shephard's question. It will
be shown that a ``random'' $n$-dimensional subspace of $\ell^{2n}_\infty$
has a unit ball, all of whose shadows are very large compared with those of
a Euclidean ball of the same volume. Such examples suggest that the
Shephard problem is less delicate than the Busemann-Petty problem for
sections: \ it is an important open question as to whether there are highly
pathological examples for the latter problem. This question is usually
referred to as the slicing problem.\vfill\eject

\noindent {\bf \S 1. A local isoperimetric inequality.}

In [V], Vaaler conjectured that for every $n\in {\Bbb N}$, every symmetric
convex body $C\subset {\Bbb R}^n$ and every $k<n$, there is an affine image
$TC$ of $C$ (for some automorphism $T$ of ${\Bbb R}^n$) so that for every
$k$-dimensional subspace $H$ of ${\Bbb R}^n$

$$|H\cap TC| \ge |TC|^{k\over n}.$$

\noindent  (Vaaler actually conjectured something slightly stronger, which is
false for small values of $k$.) This conjecture strengthens the slicing
mentioned above,
namely: \ there exists a $\delta >0$ so that for each $n$ and $C$ there
is a 1-codimensional subspace $H$ of ${\Bbb R}^n$ with

$$|H\cap C| \ge \delta|C|^{n-1\over n}.$$

The case $k=1$ of Vaaler's conjecture (for arbitrary $n$ and $C$) was
proved by the present author in [B$_1$], where the result is stated as a volume
ratio estimate. The proof
for $k=1$ really estimates volumes of 1-dimensional shadows and then uses the
fact that the smallest 1-dimensional section of a convex body is its
smallest 1-dimensional shadow.	Since minimal sections and minimal shadows
are not identical for subspaces of dimension larger than 1, such an
argument cannot be employed if $k>1\colon$ \ but it is natural to ask
whether Vaaler's conjecture can be proved for shadows (of dimension other
than 1) independently of the outstanding problem for sections. The
principal result of this paper deals with the most important case,
\hfill\break $k = n-1$.\medskip

\noindent {\bf Theorem 1.} Let $C$ be a convex body in ${\Bbb R}^n$. There
is an affine image $\widetilde C$ of $C$ (with non-zero volume) so that for
each unit vector $\theta \in {\Bbb R}^n$,

$$|P_\theta \widetilde C| \ge |\widetilde C|^{n-1\over n}.$$

\noindent The result is exactly best possible as shown by the cube.

The proof of  Theorem 1 uses the well-known theorem of John [J] which
characterises ellipsoids of minimal volume containing convex bodies. This
result is stated here as a lemma.

\noindent {\bf Lemma 2.} Let $K$ be a symmetric convex body in ${\Bbb
R}^n$. The ellipsoid of minimal volume containing $K$ is the Euclidean unit
ball $B^n_2$ if and only if $K$ is contained in $B^n_2$ and
 there are Euclidean unit vectors $(u_i)^m_1$
(for some $m\in {\Bbb N})$ on the boundary $\partial K$ of $K$ and positive
numbers $(c_i)^m_1$ so that

$$\sum^m_1 c_iu_i \otimes u_i = I_n.$$

\noindent (Here, $u_i \otimes u_i$ is the usual rank-1 orthogonal
projection onto the span of $u_i$ and $I_n$ is the identity on ${\Bbb
R}^n$.)  The identity above states that the $u_i$'s are distributed rather
like an orthonormal basis in that for each $x\in {\Bbb R}^n$,

$$|x|^2 = \sum^m_1 c_i \langle u_i,x\rangle^2.$$

\noindent The equality of the traces of the operators appearing above shows that
$\sum\limits ^m_1 c_i = n$.

Theorem 1 will be deduced from the following, which is little more than an
affine invariant reformulation.\medskip

\noindent {\bf Theorem 3.} Suppose $C$ is a convex body in ${\Bbb R}^n,
(u_i)^m_1$ a sequence of unit vectors in ${\Bbb R}^n$ and $(c_i)^m_1$ a
sequence of positive numbers for which

$$\sum^m_1 c_iu_i \otimes u_i = I_n.$$

\noindent For each $i$, let $P_i$ be the orthogonal projection $P_{u_i}$ along
$u_i$. Then

$$|C|^{n-1} \le \prod^m_1 |P_iC|^{c_i}.$$

There is an obvious relationship between Theorem 3 and Lemma 2. Theorem 3
is
closely related to an inequality of Brascamp and Lieb [Br-L] which has been
used in several places by this author, [B$_1$] and [B$_2$]. Theorem 3 and
generalisations of it were conjectured in [Br-L] (in a different form). The
special case of Theorem 3 in which the $u_i$'s form an orthonormal basis
(in which case, necessarily, $c_i = 1$ for $1\le i \le m=n$) was proved by
Loomis and Whitney [L-W]. Theorem 3 can be regarded as an isoperimetric
inequality in that it estimates the volume of a body in terms of an average
of volumes of its shadows: \ in this case, a geometric average. The key
point is that the ``number'' of shadows involved is small enough that the
local isoperimetric inequality of Theorem 1 can be deduced.

\noindent {\bf Proof of Theorem 1.}
 Because of the intertwining property of $\Pi^*$ with linear
transformations (1), there is an affine image $\widetilde C$ of $C$ so that
the ellipsoid of minimal volume containing $\Pi^*(\widetilde C)$ is the
Euclidean ball $B^n_2$. This ensures that

$$|P_\theta (\widetilde C)| \ge 1$$

\noindent for every unit vector $\theta \in {\Bbb R}^n$ and, by Lemma 2,
that there are unit vectors $(u_i)^m_1$ and positive numbers $(c_i)^m_1$ so
that

$$|P_{u_i} \widetilde C| = 1 \quad \hbox{for each}\quad i$$

\noindent and $\sum\limits ^m_1 c_iu_i \otimes u_i = I_n$.\medskip

\noindent Now, from Theorem 3,

$$|\widetilde C|^{n-1} \le \prod ^m_1 |P_{u_i} \widetilde C|^{c_i} = 1.$$

$\hfill \square$

The proof of Theorem 3 uses Minkowski's inequality for mixed volumes to
establish a duality between the 1-codimensional problem to be solved and a
1-dimensional problem. The relevant information on mixed volumes is
included here for completeness.

For a fixed $n$, let ${\cal C} = {\cal C}_n$ be the set of compact, convex
subsets of ${\Bbb R}^n$. ${\cal C}$ can be regarded as a convex cone under
Minkowski addition and multiplication by non-negative scalars. A crucial
theorem of Minkowski states that there is a symmetric $n$-positive-linear
form

$$V\colon \ \underbrace{{\cal C}\times \cdots \times{\cal C}}_{n \ {\rm
times}} \longrightarrow [0,\infty)$$

\noindent whose diagonal is volume: \ i.e. $V$ is positive linear in each
of its $n$ arguments and, for each $C\in {\cal C}$,

$$|C| = V(C,\ldots, C).$$

\noindent The values of $V$ are called mixed volumes. As a consequence of
Minkowski's theorem, the volume $|C+tD|, (t\in [0,\infty))$ can be expanded
as a polynomial in $t$,

$$\eqalignno{|C+tD| &= \sum ^n_0 {n\choose k} v_{n-k} (C,D)t^k&(2)\cr
\noalign{\hbox{where}}
v_{n-k}(C,D) &= V(\underbrace{C,\ldots, C}_{n-k}, \underbrace{D,\ldots,
D}_{k})}$$

\noindent and is called the $n-k^{\rm th}$ mixed volume of $C$ and $D$.

The Brunn-Minkowski inequality states that

$$|C+tD|^{1\over n}
$$

\noindent is a concave function of $t$ (on $[0,\infty)$). Differentiation
of (2) at $t=0$ gives Minkowski's inequality

$$|C|^{n-1\over n} |D|^{1\over n} \le v_{n-1} (C,D).\eqno (3)$$

\noindent If $D = B^n_2$ is the Euclidean unit ball, (3) is the classical
isoperimetric inequality. Inequality (3) will be used here with an
appropriate choice of $D$.

A Minkowski sum of line segments

$$\sum ^m_1[-x_i,x_i] = \bigg\{ x\in
{\Bbb R}^n\colon \ x  = \sum^m_1 \lambda_ix_i \ \hbox{ for some sequence} \
(\lambda_i)^m_1 \ {\rm with} \ |\lambda_i| \le 1, 1 \le i \le m\bigg\}$$

\noindent is
called
a zonotope. It is easily checked that if $u$ is a unit vector and $D = [-u,u]
= \{x\in {\Bbb R}^n\colon \ x = \lambda u$ for some $\lambda \in [-1,1]\}$
then

$$v_{n-1} (C,D) = {2\over n} |P_uC|$$

\noindent for any convex body $C$. So if $Z$ is the zonotope

$$Z = \sum ^m_1 \alpha_i[-u_i,u_i]$$

\noindent with $(u_i)^m_1$ a sequence of unit vectors and $(\alpha_i)^m_1$
a sequence of positive numbers,

$$v_{n-1}(C,Z) = {2\over n} \sum ^m_1 \alpha_i|P_{u_i} C|.\eqno (4)$$

\noindent There is equality in Minkowski's inequality (3) if $C=D$ and so
with $Z$ as before

$$|Z| = {2\over n} \sum ^m_1 \alpha_i |Pu_i Z|.\eqno (5)$$

\noindent This identity is usually called the volume formula for zonotopes.
A simple induction argument can be used to obtain an expression for $|Z|$
in terms of the $\alpha_i$'s and determinants of square matrices formed
from the $u_i$'s. In the following lemma, a similar inductive argument is
used to obtain an estimate for $|Z|$ which is easier to use than the actual
value.

\noindent {\bf Lemma 4.} Let $(u_i)^m_1$ be a sequence of unit vectors in
${\Bbb R}^n, (c_i)^m_1$ a sequence of positive numbers with

$$\sum ^m_1 c_iu_i \otimes u_i = I_n$$

\noindent and $(\alpha_i)^m_1$ another sequence of positive numbers. If

$$\eqalignno{Z &= \sum ^m_1 \alpha_i[-u_i, u_i]\cr
\noalign{\hbox{then}}
|Z| &\ge 2^n \prod^m_1 \Big({\alpha_i\over c_i}\Big)^{c_i}.}$$

\noindent {\bf Proof.} The proof uses induction on the dimension $n$. For
$n=1$,

$$\eqalign{|Z| &= 2 \sum ^m_1 \alpha_i\cr
&= 2 \sum^m_1 c_i\Big({\alpha_i\over c_i}\Big)\cr
&\ge 2 \prod ^m_1 \Big({\alpha_i\over c_i}\Big)^{c_i}}$$

\noindent by the AM-GM inequality since $\sum\limits ^m_1 c_i = 1$, when
$n=1$.

For larger $n$, the volume formula (5) shows that with $P_i = P_{u_i}$

$$\eqalignno{|Z| &= {2\over n} \sum ^m_1 \alpha_i |P_iZ|\cr
&= 2 \sum ^m_1 {c_i\over n} {\alpha_i\over c_i} |P_iZ|\cr
&\ge 2 \prod ^m_1 \Big({\alpha_i\over c_i} |P_iZ|\Big)^{c_i\over n}&(6)}$$

\noindent since $\sum\limits ^m_1 c_i = n$.

Now, for each fixed $i, P_iZ$ is a zonotope with summands

$$\alpha_j[-P_i(u_j), P_i(u_j)],\qquad 1 \le j \le m$$

\noindent and contained in the $(n-1)$-dimensional space $P_i({\Bbb R}^n)$.
For each $i$ and $j$ let

$$\gamma_{ij} = |P_i(u_j)| = |P_j(u_i)|$$

\noindent (so $\gamma^2_{ij} = 1-\langle u_i,u_j\rangle^2$).  Then for each
$i$,

$$P_iZ = \sum^m_{j=1} \alpha_j \gamma_{ij} [-v_{ij}, v_{ij}]$$

\noindent where $v_{ij}$ is the unit vector in the direction of $P_iu_j$
(or any direction if $\gamma_{ij} = 0$). Now, for each $i$,

$$\eqalign{P_i &= \sum^m_{j=1} c_jP_iu_j \otimes P_iu_j\cr
&= \sum ^m_{j=1} c_j \gamma_{ij}^2 v_{ij} \otimes  v_{ij}}$$

\noindent and $P_i$ acts as the identity on $P_i({\Bbb R}^n)$. So, by the
inductive hypothesis,

$$|P_iZ| \ge 2^{n-1} \prod ^m_{j=1} \Big({\alpha_j\gamma_{ij}\over c_j
\gamma^2_{ij}}\Big)^{c_j\gamma^2_{ij}}$$

\noindent where it is understood that if $\gamma_{ij} = 0$, the $j^{\rm
th}$ factor is 1:  \ (so in particular, the $i^{\rm th}$ factor is 1).

Substitution of the inequalities for each $i$ into (6) shows that

$$|Z| \ge 2^n \bigg(\prod^m_{i,j=1} \Big({\alpha_i\over c_i}\Big)^{c_i}
\Big({\alpha_j\over
c_j\gamma_{ij}}\Big)^{c_ic_j\gamma^2_{ij}}\bigg)^{1\over n}$$

\noindent and this expression is at least

$$2^n \prod ^m_{i=1} \Big({\alpha_i\over c_i}\Big)^{c_i}$$

\noindent because ${1\over \gamma_{ij}}\ge 1$ for all $i$ and $j$ and, for
each $j$,

$$\sum ^m_{i=1} c_i\gamma^2_{ij} = \sum ^m_{i=1} c_i(1 - \langle
u_i,u_i\rangle^2) = n-1.$$

$\hfill \square$

\noindent {\bf Proof of Theorem 3.} Let $C, (u_i)^m_1, (c_i)^m_1$ and
$(P_i)^m_1$ be as in the theorem's statement. For each $i$, set

$$\alpha_i = {c_i\over |P_iC|}.$$

\noindent With $Z = \sum\limits ^m_1 \alpha_i[-u_i,u_i]$, (3) and (4)
above, show that

$$\eqalign{|C|^{n-1} &\le |Z|^{-1} \bigg( {2\over n} \sum ^m_1 \alpha_i
|P_iC|\bigg)^n\cr
&= |Z|^{-1} \bigg({2\over n} \sum ^m_1 c_i\bigg)^n\cr
&= 2^n |Z|^{-1}}$$

\noindent since $\sum\limits ^m_1 c_i = n$.\medskip

\noindent Now by Lemma 4,

$$|C|^{n-1} \le 2^n \bigg(2^n \prod ^m_{i=1} \Big({\alpha_i\over
c_i}\Big)^{c_i}\bigg)^{-1} = \prod ^m_{i=1} |P_iC|^{c_i}.\eqno \square$$

\noindent {\bf Remark 1.}  The lower estimate for volumes of zonotopes
given by Lemma 4 could have been obtained from the volume ratio estimates
proved in the author's paper [B$_2$] together with Reisner's reverse
Santalo inequality for zonoids, [R].

\noindent {\bf Remark 2.} For every  convex body $C$ in ${\Bbb R}^n$, there
is an affine image $\widetilde C$ of $C$ with

$$|P_\theta \widetilde C| \le M\sqrt{n} |\widetilde C|^{n-1\over n}$$

\noindent for each unit vector $\theta$ ($M$ being an absolute constant).
This estimate depends upon the fact that subspaces of $L_1$ have uniformly
bounded volume ratios. It seems likely that $M$ could be taken to be 1 for
symmetric convex bodies $C$ (the cube again being extremal). The estimates
in [B$_2$] show that $M$ can be taken to be ${2\sqrt{e}\over \pi} \approx
1.05$ in this case.\vfill\eject

\noindent {\bf \S 2. A remark on the Shephard problem.}

This section contains a strongly negative solution to the Shephard problem
described in the introduction. Petty and Schneider ([P] and [S])
constructed pairs of bodies $C$ and $D$ in ${\Bbb R}^n$ so that

$$|P_\theta D| \le |P_\theta C| \quad \hbox{for all} \quad \theta \in
S^{n-1},\eqno (7)$$

\noindent  but $|D| > |C|$. Schneider also showed that the conclusion $|D|
\le |C|$ does hold if $C$ is a zonoid (a limit, in the Hausdorff metric, of
zonotopes). A little more generally, if $C$ and $D$ are convex bodies
satisfying (7) and
$Z = \sum\limits^m_1 \alpha_i[-u_i,u_i]$ is a zonotope included in $C$
(with $u_i \in S^{n-1}, 1 \le i\le m$), then, rather as in the proof of
Theorem 3,

$$\eqalign{|D|^{n-1\over n} |Z|^{1\over n} &\le v_{n-1} (D,Z) = \sum^m_1
\alpha_i|P_{u_i} D|\cr
&\le \sum^m_1 \alpha_i |P_{u_i} C| = v_{n-1}(C,Z)\cr
&\le v_{n-1} (C,C) = |C|}$$

\noindent where the last inequality is a consequence of the monotonicity of
mixed volumes: \ (this particular case is obvious from the fact that if $Z
\subset C$,

$$C + tZ \subset (1+t)C$$

\noindent for all $t\ge 0$).

Hence,

$$|D| \le \Big({|C|\over |Z|}\Big)^{1\over n-1} |C|.$$

\noindent  Since every convex body $C$ contains a zonoid $Z$ with
${|C|\over |Z|}$ not too large, the last inequality shows that under
hypothesis (7), one does have

$$|D| \le {3\over 2} \sqrt{n} |C|.$$

\noindent Apart from a constant factor, it turns out that this is the most
that can be said. The Euclidean ball of volume 1 in ${\Bbb R}^n$ has
shadows
of 1-codimensional volume about $\sqrt{e}$ (as $n\to \infty$): \ in Theorem
5 it is shown that there is a body of volume 1 in ${\Bbb R}^n$, all of
whose shadows have volume about $\sqrt{n}$.\medskip

\noindent {\bf Theorem 5.} There is a constant $\delta >0$ so that for each
$n \in {\Bbb N}$, there is a symmetric convex body $K$ in ${\Bbb R}^n$
satisfying,

$$|P_\theta K| \ge \delta \sqrt{n} |K|^{n-1\over n}$$

\noindent for every unit vector $\theta \in {\Bbb R}^n$.

The proof of Theorem 5 depends heavily upon the theorem of Minkowski on the
existence of bodies with given projections. An appropriately detailed
statement of this theorem is given as Lemma 6 below: \ the following
notation is needed. For a sequence $(u_i)^m_1$ of unit vectors spanning
${\Bbb R}^n$ and a sequence $(\gamma_i)^m_1$ of positive numbers let
${\cal F} = {\cal F}((u_i), (\gamma_i))$ be the family of convex bodies of
the form

$$\{x\in {\Bbb R}^n\colon \ |\langle x,u_i\rangle| \le t_i, \quad 1 \le i
\le m\}$$

\noindent indexed by sequences $(t_i)^m_1$ of positive reals satisfying

$$\sum ^m_1 \gamma_it_i = 1.$$

\noindent {\bf Lemma 6.} With the above notation, ${\cal F}$ has an unique
element of maximal volume, $K$ (say), satisfying

$$|P_\theta K| = {n|K|\over 2} \sum ^m_1 \gamma_i |\langle u_i,\theta
\rangle |$$

\noindent for each $\theta \in S^{n-1}$.$\hfill \square$

The body that satisfies the conclusion of Theorem 5 will be the unit ball
of a ``random'' $n$-dimensional subspace of $\ell^{2n}_\infty$. Such
subspaces are known to have many pathological properties stemming from the
fact that the $\ell ^{2n}_1$ and $\ell ^{2n}_2$ norms are well-equivalent
on such spaces: \ this was proved in [F-L-M]. For the history of this
result and its many extensions, see e.g. [M-S]. The form of the result
needed here is given as a Lemma.

\noindent {\bf Lemma 7.} There is a $\delta >0$ so that for $n\in {\Bbb N}$
there are unit vectors $u_1,\ldots, u_{2n}$ in ${\Bbb R}^n$ with

$$\sum ^{2n}_1 |\langle x,u_i\rangle | \ge \delta\sqrt{n}|x|$$

\noindent for every vector $x\in {\Bbb R}^n$.$\hfill \square$

The last lemma that will be needed is a result of Vaaler [V], concerning
volumes of sections of the cube in ${\Bbb R}^n$. The form required here is
the following.

\noindent  {\bf Lemma 8.} Let $(u_i)^m_1$ be a sequence of unit vectors in
${\Bbb R}^n$. Then the volume of the symmetric convex body with these
vectors as boundary functionals satisfies.

$$|\{x\in {\Bbb R}^n\colon \ |\langle x,u_i\rangle|\le 1, \quad 1 \le i \le
m\}|^{1\over n} \ge 2\sqrt{n\over m}.\eqno \square$$

\noindent {\bf Proof of Theorem 5.} Let $(u_i)^{2n}_1$ be a sequence of
vectors with the property described in Lemma 7. Take $m = 2n$ and apply
Lemma 6 with $\gamma_i = {1\over m}, 1 \le i \le m$.

The body $K$ of maximal volume in the family ${\cal F}$ satisfies

$$|P_\theta K| = {n|K|\over 2m} \sum ^m_1 |\langle u_i,\theta\rangle| \ge
{|K|\over 4} \cdot \delta \sqrt{n}$$

\noindent for every unit vector $\theta$. Now, ${\cal F}$ also contains
the body

$$C = \{x\in {\Bbb R}^n\colon \ |\langle x,u_i\rangle| \le 1,\quad 1 \le i
\le m\}.$$

\noindent So (by the maximality of $K$), $|K| \ge |C|$. But, by Lemma 8,

$$|C|^{1\over n} \ge 2\sqrt{n\over m} = \sqrt{2}$$

\noindent and hence $|K|^{1\over n} \ge \sqrt{2}$. Therefore, for every
unit vector $\theta$,

$$|P_\theta K| \ge {\delta \sqrt{n}\over 2\sqrt{2}} |K|^{n-1\over n}.\eqno
\square$$

\noindent {\bf Remark.}  The above argument can be extended slightly to
estimate surface area to volume ratios of subspaces of $\ell_\infty$. If
$iq(X)$ is the isoperimetric quotient of the finite-dimensional normed
space $X$ (see [Sch\"u] for definitions) normalised so that $iq(\ell^n_2)
= 1$ then for every $n$-dimensional subspace $X$ of $\ell^m_\infty$,

$$iq(X) \ge {\delta \sqrt{n}\over \sqrt{1+\log {m\over n}}}$$

\noindent for some absolute constant  $\delta >0$.\vfill\eject

\noindent {\bf References.}

\item{[B$_1$]} K.M. Ball, Volumes of sections of cubes and related
problems, Israel Seminar \hfill\break (G.A.F.A.) 1988, Springer-Verlag,
Lecture Notes \#1376, (1989), 251-260.

\item{[B$_2$]} K.M. Ball, Volume ratios and a reverse isoperimetric
inequality, In preparation.

\item{[B-L]} J. Bourgain and J. Lindenstrauss, Projection bodies, Israel
Seminar (G.A.F.A.) 1986-87, Springer-Verlag, Lecture Notes \#1317, (1988),
250-269.

\item{[Br-L]} Herm Jan Brascamp and Elliott H. Lieb, Best constants in
Young's inequality, its converse and its generalization to more than three
functions, Advances in Math. 20 (1976), 151-173.

\item{[B-P]} H. Busemann and C.M. Petty, Problems on convex bodies, Math.
Scand. 4 (1956), 88-94.

\item{[F-L-M]} T. Figiel, J. Lindenstrauss and V.D. Milman, The dimension
of almost spherical sections of convex bodies, Acta Math. 139 (1977),
53-94.

\item{[J]} F. John, Extremum problems with inequalities as subsidiary
conditions, Courant Anniversary Volume, Interscience, New York, 1948,
187-204.

\item{[L-R]} D.G. Larman and C.A. Rogers, The existence of a centrally
symmetric convex body with central sections that are unexpectedly small,
Mathematika 22 (1976), 164-175.

\item{[L-W]} L.H. Loomis and H. Whitney, An inequality related to the
isoperimetric inequality, Bull. Amer. Math. Soc. 55 (1949), 961-962.

\item{[M-S]} V.D. Milman and G. Schechtman, Asymptotic Theory of Finite
Dimensional Normed Spaces, Springer-Verlag, Lecture Notes \#1200, (1986),
64-104.

\item{[P]} C.M. Petty, Projection bodies, Proc. Colloq. on Convexity,
Copenhagen (1967), 234-241.

\item{[R]} S. Reisner, Zonoids with minimal volume product, Math. Z. 192
(1986), 339-346.

\item{[S]} R. Schneider, Zu einem Problem von Shephard \"uber die
Projektionen konvexer \hfill\break K\"orper, Math. Z. 101 (1967), 71-82.

\item{[Sch\"u]} C. Sch\"utt, The isoperimetric quotient and some classical
Banach spaces, to appear.

\item{[V]} J.D. Vaaler, A geometric inequality with applications to linear
forms, Pacific J. Math. 83 (1979), 543-553.

\end